\documentclass[14pt]{article}
\usepackage{amsfonts}
\usepackage{latexsym}
\usepackage{cite}
\usepackage{amsmath,amsfonts,latexsym,amssymb}
\usepackage[mathscr]{eucal}
\usepackage{cases,color}
\usepackage{amsthm}

\usepackage[bf,small]{caption2}
\usepackage{float}
\usepackage{graphicx}
\usepackage{amsmath}
\usepackage{amssymb}
\usepackage[all]{xy}

\newtheorem{theorem}{theorem}[section]
\newtheorem{thm}[theorem]{Theorem}
\newtheorem{lem}[theorem]{Lemma}
\newtheorem{prob}[theorem]{Problem}

\newtheorem{nota}[theorem]{Notation}
\newtheorem{defn}[theorem]{Definition}

\newtheorem{rmk}[theorem]{Remark}

\begin{document}

\title{\textbf{Computing untwisted Dijkgraaf-Witten invariants for arborescent links}}
\author{\Large Haimiao Chen}

\date{}
\maketitle

\begin{abstract}
  We briefly review 3-dimensional untwisted Dijkgraaf-Witten theory over a finite group $\Gamma$, and present a method of computing untwisted Dijkgraaf-Witten invariants for arborescent links. Some explicit formulas are given when $\Gamma=\mathbb{Z}/p\mathbb{Z}\rtimes\mathbb{Z}/(p-1)\mathbb{Z}$ for an odd prime $p$.

\end{abstract}

\section{Introduction}

Dijkgraaf-Witten theory was first proposed in \cite{DW90}, and since then has been
further studied by many people; see \cite{Fe93, FQ93, Fr94-1, Fr94-2, Mo15, Wa92, Ye92} and the references therein.
It can be defined for any integer $d\geqslant 1$, any finite group $\Gamma$ and any cohomology class
$\alpha\in H^{d}(\Gamma;U(1))$, and is called untwisted when $\alpha=0$. The key ingredient of the untwisted DW theory
(denoted as ${\rm DW}^{0}_{d,\Gamma}$) is counting homomorphisms from the fundamental group of a manifold to $\Gamma$.

DW theory not only has theoretical importance since it is one of the first rigorously constructed TQFTs, but also is practically
interesting because of its relation with fundamental group. The fundamental group of a manifold $M$, as a noncommutative object, encodes
much topological information of $M$, but at the same time, is often hard to handle. The DW invariant of $M$,
however, can extract partial information on $\pi_1(M)$.

In general, enumeration of homomorphisms from a group $G$ to finite groups is usually helpful for understanding $G$.
For example, the classical topic of counting finite-index subgroups of a finitely generated group $G$ is
directly related to counting homomorphisms from $G$ to symmetric groups, see  \cite{Lu95, Su01}; the Hall
invariant, which counts epimorphisms from $G$ onto a finite group, dates back to \cite{Ha36} and also has its own interests
\cite{MS02, Su01}.

When $G=\pi_1(M)$, homomorphisms from $G$ to finite groups have additional topological meanings.
In dimension 2, people studied existence and enumeration problems on surface coverings by counting homomorphisms from surface groups
to finite groups, see \cite{Jo95, PP06, Tu09}. But in 3-dimension, till now, besides \cite{LM00} and the author's work \cite{Ch12-1}, there are few such results. In knot theory, for a knot $K\subset S^3$, there are good reasons for paying attention to homomorphisms from
$\pi(K):=\pi_{1}(E_K)$ (where $E_K$ is the complement of a tubular neighborhood of $K$) to finite groups. It was shown in \cite{Ei07} that all finite quandle invariants can be expressed in terms of the numbers $N_{\Gamma}(x,y)$ of homomorphisms $\pi(K)\rightarrow\Gamma$ sending the
meridian and the longitude to $x$ and $y$, respectively, for appropriate finite group
$\Gamma$ and various pairs $(x,y)\in\mathfrak{P}_\Gamma:=\{(x,y)\colon xy=yx\}$; it was also pointed out that such enumeration is helpful to somehow
keep track of the knot group since the knot group is residually finite. In fact the DW invariant of $E_K$ is a specific vector comprising exactly the numbers $N_{\Gamma}(x,y)$.

In recent years, many deep connections between 3-manifolds and finite groups have been revealed (see \cite{BF20,Ko13,LR10,Wi17,WZ17} and the related references). There are a lot of interesting things to be further explored.
In \cite{Ch18}, the author deduced a criterion for a link $L$ to be periodic, in terms of the DW invariants of $L$ as well as that of the quotient link. For the criterion to be applicable, we need some general knowledge on DW invariants of the ``nonvisible" quotient link.
It is thus natural to formulate the following.
\begin{prob}
Given a finite group $\Gamma$, which function
$$\mathfrak{P}_\Gamma^n=\{(x_1,h_1;\ldots;x_n,h_n)\colon x_ih_i=h_ix_i,\ 1\le i\le n\}\to\mathbb{N}$$
can be realized as the Dijkgraaf-Witten invariant of some $n$-component link?
\end{prob}

In the study of DW theory, a big problem is the difficulty of practical and concrete computations.

As a widely known fact \cite{RT91,Ba01}, in general each {\it modular tensor category} $\mathcal{C}$ gives rise to a $3$-dimensional extended TQFT ${\rm RT}_{\mathcal{C}}$. When $\mathcal{C}$ is $\mathcal{E}(\Gamma)$ (to be introduced in Section \ref{sec:cat}), a result of Freed \cite{Fr94-2} established the equivalence between ${\rm DW}^{0}_{3,\Gamma}$ and the Reshetikhin-Turaev theory ${\rm RT}_{\mathcal{E}(\Gamma)}$. So the DW invariants of links can be computed via diagrams as for general RT invariants. However, there is no such computation seen in the publications.

The main contribution of the present paper consists of two parts. First, based on the preprint \cite{Ch12-2}, we clarify several aspects of DW theory which were usually folklore or spreading among experts, and present a practical method for computing DW invariants for arborescent links. Second, we provide significant ingredients of the MTC $\mathcal{E}(\Gamma)$ when $\Gamma$ is the semi-direct product $\mathbb{Z}/p\mathbb{Z}\rtimes\mathbb{Z}/(p-1)\mathbb{Z}$, where $p$ is an odd prime number; we will see that much number-theoretic richness is reflected in the structure of $\mathcal{E}(\Gamma)$.

The content is organized as follows.
Section 2 is a preliminary on DW theory, tangle and link, and RT theory.
In Section 3 we develop a general method of computing DW invariants for arborescent links.
In Section 4, we give key formulas for $\Gamma=\mathbb{Z}/p\mathbb{Z}\rtimes\mathbb{Z}/(p-1)\mathbb{Z}$, and as an illustration, we compute the DW invariants for a family of knots.

\medskip

\begin{nota}
\rm
For positive integers $m,n$, let $\mathcal{M}(m,n)$ denotes the space of $m\times n$ matrices over $\mathbb{C}$.

For a set $X$, let $\#X$ denote its cardinality.

Let $\mathbb{T}$ denote the torus $S^{1}\times S^{1}$. Let $\mathfrak{m}_0=S^{1}\times 1$ and $\mathfrak{l}_0=1\times S^{1}$.

For a group $G$, usually denote the identity element by $e$. For $x,y\in G$, let $x\lrcorner y$ denote $xyx^{-1}$;
let $\delta^x_y=1$ if $x=y$ and $\delta^x_y=0$ otherwise;
let ${\rm Cen}(x)$ denote the centralizer of $x$; let ${\rm Con}(x)$ denote the conjugacy class containing $x$.

For a condition $\mathfrak{x}$, let $\delta_{\mathfrak{x}}=1$ if $\mathfrak{x}$ holds, and $\delta_{\mathfrak{x}}=0$ otherwise.

If $x\in G$ has order $n$, and $k\in\mathbb{Z}/n\mathbb{Z}$, then by $x^k$ we mean $x^{\tilde{k}}$ for any $\tilde{k}\in\mathbb{Z}$ whose residue class modulo $n$ is $k$.
\end{nota}

\section{Preliminary}

\subsection{Untwisted Dijkgraaf-Witten theory} \label{sec:DW}

For a finite group and a positive integer $d$, the $d$-dimensional untwisted Dijkgraaf-Witten theory $Z={\rm DW}^0_{d,\Gamma}$ over $\Gamma$ is given as follows.

For a connected closed $(d-1)$-manifold $B$, let $Z(B)$ be the vector space of functions $\alpha:\hom(\pi_{1}(B),\Gamma)\to\mathbb{C}$ such that $\alpha(h\lrcorner\phi)=\alpha(\phi)$ for all $h\in\Gamma$ and $\phi\in \hom(\pi_{1}(B),\Gamma)$, where $(h\lrcorner\phi)(x)=h\lrcorner\phi(x)$.
For connected closed $(d-1)$-manifolds $B_1,\ldots,B_r$, set $Z(\sqcup_{i=1}^rB_i)=\bigotimes_{i=1}^rZ(B_i)$.

For a connected closed $d$-manifold $C$, let
\begin{align*}
Z(C)=\frac{1}{\#\Gamma}\cdot\#\hom(\pi_{1}(C),\Gamma)\in\mathbb{C}=:Z(\emptyset);
\end{align*}
for a connected $d$-manifold $C$ with $\partial C\ne\emptyset$, define $Z(C)\in Z(\partial C)$ by setting
\begin{align*}
Z(C)(\phi)=\#\{\Phi\in\hom(\pi_{1}(C),\Gamma)\colon\Phi|_{\partial C}=\phi\};
\end{align*}
for connected $d$-manifolds $C_1,\ldots,C_s$, set
$$Z(\sqcup_{i=1}^sC_i)=\otimes_{i=1}^sZ(C_i)\in{\bigotimes}_{i=1}^sZ(\partial C_i).$$

When $d=3$, the vector space $E:=Z(\mathbb{T})$ can be identified with the space of functions
$\alpha:\mathfrak{P}_\Gamma\to\mathbb{C}$ such that $\alpha(a\lrcorner x,a\lrcorner g)=\alpha(x,g)$ for all $a\in\Gamma$.

For each conjugacy class $\mathfrak{c}$ of $\Gamma$, choose $x_{\mathfrak{c}}\in \mathfrak{c}$, take an irreducible representation $\rho$ of the centralizer ${\rm Cen}(x_{\mathfrak{c}})$ of $x_{\mathfrak{c}}$, and put
$$\chi_{\mathfrak{c},\rho}(y,b)=\begin{cases} {\rm tr}(\rho(g^{-1}\lrcorner b)), &y=g\lrcorner x_{\mathfrak{c}}, \\ 0,&y\notin\mathfrak{c}. \end{cases}$$
By Lemma 5.4 of \cite{Fr94-1}, all such $\chi_{\mathfrak{c},\rho}$'s form a canonical orthonormal basis for $E$,
with respect to the inner product
$$(\alpha_1,\alpha_2)_E=\frac{1}{\#\Gamma}\cdot\sum_{(x,g)\in\mathfrak{P}_\Gamma}\alpha_1(x,g)\overline{\alpha_2(x,g)}.$$

\subsection{Framed Tangles and links} \label{section:tangle}

A {\it framed tangle} $T$ is an equivalence class of pairs $(\underline{T},{\rm fr})$ where $\underline{T}$ is an oriented 1-manifold embedded in
$[0,1]\times\mathbb{R}^{2}\subset\mathbb{R}^{3}$
such that $\partial \underline{T}\subset\{0,1\}\times\mathbb{R}\times\{0\}$, and ${\rm fr}$ is a nonzero section of the normal vector field on $\underline{T}$. Two pairs $(\underline{T},{\rm fr})$ and $(\underline{T}',{\rm fr}')$ are equivalent if there is an orientation-preserving homeomorphism of $\mathbb{R}^{3}$ that takes $\underline{T}$ to $\underline{T}'$ and takes ${\rm fr}$ to ${\rm fr}'$.
Denote $T=[\underline{T},{\rm fr}]$. We do not distinguish a framed tangle from its representatives. Call $s(T):=(\{0\}\times\mathbb{R}\times\{0\})\cap\underline{T}$ (resp. $t(T):=(\{1\}\times\mathbb{R}\times\{0\})\cap \underline{T}$) the
{\it source} (resp. the {\it target}) of $T$. Regard $s(T)$ and $t(T)$ as oriented 0-manifolds.

A tangle with empty source and target is just an oriented link.

\begin{figure}[h]
  \centering 
  \includegraphics[width=4.6cm]{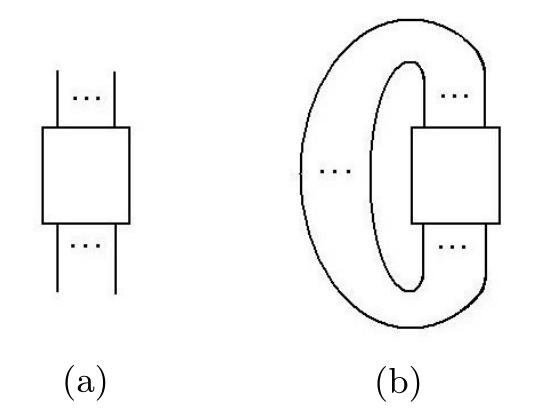}\\
  \caption{(a) A tangle $T$; (b) the closure of $T$.} \label{figure:closure}
\end{figure}

\begin{figure}[h]
  \centering 
  \includegraphics[width=6cm]{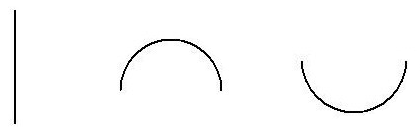}\\
  \caption{Three of the basic tangles.} \label{figure:basic}
\end{figure}

\begin{figure}[h]
  \centering 
  \includegraphics[width=9cm]{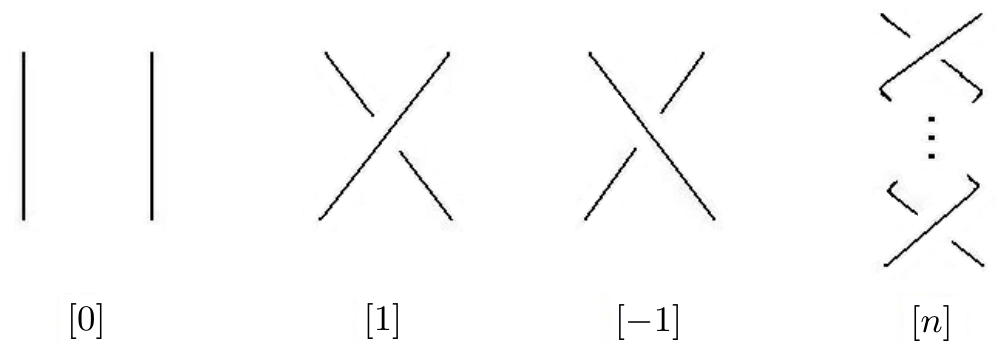}\\
  \caption{Integral tangles.} \label{figure:tangle}
\end{figure}

A framed tangle can always be presented as a tangle diagram, with the blackboard framing understood. In this vein, two tangle diagrams
represent the same framed tangle if and only if their underlying tangles are equivalent and the writhe numbers of the corresponding components are all equal.

Given tangles $T_{1},T_{2}$, the {\it horizontal composite} $T_{1}+T_{2}$ is defined by putting $T_{2}$ on the right
side of $T_{1}$. When $\#s(T_{2})=\#t(T_{1})$ (so that there is an obvious bijection between $s(T_{2})$ and $t(T_{1})$), and the orientations
are compatible, we can define the {\it vertical composite} $T_{2}\ast T_{1}$ to be framed tangle obtained by
putting $T_2$ on top of $T_1$ and identifying $t(T_{1})$ with $s(T_{2})$; clearly, $s(T_{2}\ast T_{1})=s(T_{1})$ and $t(T_{2}\ast T_{1})=t(T_{2})$.

When $s(T)$ can be identified with $t(T)$, the {\it closure} of $T$, denoted $\overline{T}$, is the closed 1-manifold obtained by connecting the source and target using parallel lines in the way shown in Figure \ref{figure:closure}.

Call the tangles in Figure \ref{figure:basic} and the $[\pm 1]$ shown in Figure \ref{figure:tangle} {\it basic tangles}.

For $n\in\mathbb{Z}$, the {\it integral tangle} $[n]$ is the vertical composite of $|n|$ copies of $[1]$ or $[-1]$, depending on the sign of $n$.

\begin{defn}  \label{defn:DW-of-link}
\rm Let $L=K_1\sqcup\cdots\sqcup K_n$ be a $n$-component framed oriented link. Let $\mathcal{N}(L)$ be a tubular neighborhood.
For each $i$, let $\mathfrak{l}_i=\acute{K_i}$, the knot obtained by moving $K_i$ a short distance according to the framing, and choose a meridian $\mathfrak{m}_i$ paired with $\mathfrak{l}_i$.
The {\it Dijkgraaf-Witten invariant} $Z(L)$ is the image of $Z(S^{3}\setminus\mathcal{N}(L))$ under the isomorphism $Z(\partial\mathcal{N}(L))\cong E^{\otimes n}$ induced by the homeomorphism $\partial\mathcal{N}(L)\rightarrow\sqcup^n\mathbb{T}$ sending $\mathfrak{m}_i$, $\mathfrak{l}_i$ respectively to $\mathfrak{m}_0$, $\mathfrak{l}_0$.
\end{defn}

\subsection{The category $\mathcal{E}(\Gamma)$} \label{sec:cat}

Fix a finite group $\Gamma$.
Let $\mathcal{E}(\Gamma)$ be the category of (finite-dimensional) graded vector spaces $\mathbf{u}=\bigoplus_{x\in\Gamma}\mathbf{u}_{x}$
together with a left $\Gamma$ action such that $g(\mathbf{u}_{x})=\mathbf{u}_{g\lrcorner x}$; a morphism $\mathbf{u}\rightarrow\mathbf{v}$ is a family of linear maps $f=\{f_{x}:\mathbf{u}_{x}\rightarrow \mathbf{v}_{x}\colon x\in\Gamma\}$ such that $g\circ f_{x}=f_{g\lrcorner x}\circ g$ for all $x,g$. The category $\mathcal{E}(\Gamma)$ is just the category of modules over the quantum double of $\Gamma$ \cite{Ba01}, and also the same as the category of $\Gamma$-crossed modules \cite{FL21}.

Suppose $\mathbf{u}\in\mathcal{E}(\Gamma)$, and $(x,g)\in\mathfrak{P}_\Gamma$, i.e. $xg=gx$, let
\begin{align}
\chi_{\mathbf{u}}(x,g)=\textrm{tr}(g:\mathbf{u}_{x}\rightarrow \mathbf{u}_{x}).
\end{align}
The function $\chi_{\mathbf{u}}$ belongs to $E$ and is called the {\it character} of $\mathbf{u}$.
Define
\begin{align}
\dim\mathbf{u}=\sum\limits_{x\in\Gamma}\chi_{\mathbf{u}}(x,e).  \label{eq:def-dim}
\end{align}

Let $\Lambda$ be a complete set of representatives of isomorphism classes of simple objects.
Then actually $\{\chi_{\mathbf{u}}\colon \mathbf{u}\in\Lambda\}$ is just the canonical basis $\{\chi_{\mathfrak{c},\rho}\}$ (introduced in Section \ref{sec:DW}) for $E$.
Here are the crucial structures of $\mathcal{E}(\Gamma)$:
\begin{itemize}
  \item ({\it tensor product}) A bifunctor $\odot:\mathcal{E}(\Gamma)\times\mathcal{E}(\Gamma)\rightarrow\mathcal{E}(\Gamma)$ defined by
        $$(\mathbf{u}\odot \mathbf{v})_{x}=\bigoplus\limits_{x_{1}x_{2}=x}\mathbf{u}_{x_{1}}\otimes\mathbf{v}_{x_{2}}.$$
        A ``unit" object $\mathbf{1}$ is given by $\mathbb{C}$ sitting at $e$, with $\Gamma$ acted trivially.
  \item ({\it associator}) A natural isomorphism $(\mathbf{u}\odot\mathbf{v})\odot\mathbf{w}\rightarrow \mathbf{u}\odot(\mathbf{v}\odot\mathbf{w})$
        given by the natural isomorphism of vector spaces
        $$(\mathbf{u}_{x}\otimes \mathbf{v}_{y})\otimes\mathbf{w}_{z}\cong\mathbf{u}_{x}\otimes(\mathbf{v}_y\otimes\mathbf{w}_{z}).$$
  \item ({\it braiding}) A natural isomorphism $R_{\mathbf{u},\mathbf{v}}:\mathbf{u}\odot\mathbf{v}\rightarrow \mathbf{v}\odot\mathbf{u}$ given by
        \begin{align*}
        u_{x}\otimes v_{y}\mapsto xv_{y}\otimes u_{x}, \qquad \text{for\ } u_{x}\in\mathbf{u}_{x},\ v_{y}\in\mathbf{v}_{y}.
        \end{align*}
  \item ({\it dual}) An involution $^{\star}:\mathcal{E}(\Gamma)\rightarrow\mathcal{E}(\Gamma)^{{\rm op}}$, defined by
        $$(\mathbf{u}^{\star})_{x}=(\mathbf{u}_{x^{-1}})^\ast;$$
        the action $g:\mathbf{u}^{\star}_{x}\rightarrow \mathbf{u}^{\star}_{g\lrcorner x}$ is given by the dual of $g^{-1}:\mathbf{u}_{g\lrcorner x^{-1}}\rightarrow \mathbf{u}_{x^{-1}}$.
  \item A pair of natural transformations $\iota_{\mathbf{u}}:\mathbf{1}\rightarrow \mathbf{u}\odot\mathbf{u}^{\star}$, $\epsilon_{\mathbf{u}}:\mathbf{u}^{\star}\odot \mathbf{u}\rightarrow\mathbf{1}$, which are expressed explicitly as
        \begin{align*}
        \iota_{\mathbf{u}}(1)=\sum\limits_{x}\sum\limits_{a}u_{x,a}\otimes u_{x,a}^\ast,  \qquad
        \epsilon_{\mathbf{u}}(u^\ast_{x,a}\otimes u_{y,b})=\delta^x_y\delta^a_b,
        \end{align*}
        where $\{u_{x,a}\}$ is an arbitrary basis for $\mathbf{u}_{x}$, and $\{u_{x,a}^\ast\}$ is the corresponding dual basis for $\mathbf{u}^\star_{x^{-1}}=(\mathbf{u}_x)^\ast$.
\end{itemize}

\subsection{Reshetikhin-Turaev invariant}

Suppose $\mathcal{C}$ be a modular tensor category (see \cite{Ba01} for details). It has structures which are generalizations of those possed by $\mathcal{E}(\Gamma)$ in the previous subsection. We use the same notations as for $\mathcal{E}(\Gamma)$ in Section \ref{sec:cat}; for instance, let $\odot$ denote the tensor product and let $\mathbf{1}$ denote the unit object. Let $O(\mathcal{C})$ denote the set of objects in $\mathcal{C}$. A {\it $\mathcal{C}$-colored framed tangle} is a pair $(T,c)$ consisting of a framed tangle $T$ and a map $c:\pi_{0}(T)\rightarrow O(\mathcal{C})$.

Two colored framed tangles $(T,c)$ and $(T',c')$ are considered as the same if $T'$ is obtained by reversing the orientations
of some components $A_{i},i\in I$, of $T$ and $c'$ differs from $c$ only at these components, with $c'(A_{i})\cong c(A_{i})^{\star}, i\in I$.

\begin{figure} [h]
  \centering 
  \includegraphics[width=9cm]{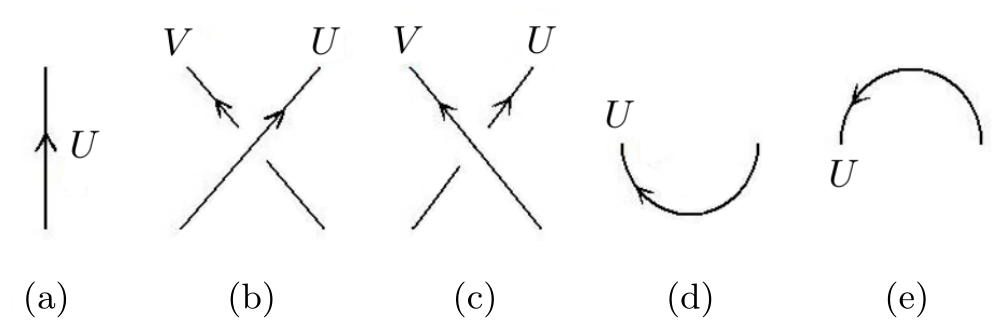}\\
  \caption{Colored basic tangles.}\label{figure:mor}
\end{figure}

Given a colored framed tangle $(T,c)$. If $s(T)=\emptyset$ (resp. $t(T)=\emptyset$), then put $\mathbf{c}_{T}=\mathbf{1}$ (resp. $\mathbf{c}^{T}=\mathbf{1}$). In general, suppose $s(T)=\{\mathfrak{s}_1,\ldots,\mathfrak{s}_k\}$, $t(T)=\{\mathfrak{t}_1,\ldots,\mathfrak{t}_\ell\}$.
For an oriented arc $A$, let $s(A)$ (resp. $t(A)$) denote its starting (resp. terminal) point.
Define $\mathbf{c}_{(i)}$ to be $c(A)$ (resp. $c(A)^{\star}$) if $\{\mathfrak{s}_{i}\}=s(A)$ (resp. $\{\mathfrak{s}_{i}\}=t(A)$) for some component $A$ of $T$, and define $\mathbf{c}^{(j)}$ to be $c(A)$ (resp. $c(A)^{\star}$) if $\{\mathfrak{t}_j\}=t(A)$ (resp. $\{\mathfrak{t}_j\}=s(A)$) for some component $A$ of $T$.
Put
\begin{align*}
\mathbf{c}_{T}=\mathbf{c}_{(1)}\odot\cdots\odot\mathbf{c}_{(k)}, \qquad
\mathbf{c}^T=\mathbf{c}^{(1)}\odot\cdots\odot\mathbf{c}^{(\ell)}.
\end{align*}
As shown in \cite{Ba01,RT91}, one can associate to $(T,c)$ a morphism
\begin{align}
\mathcal{F}(T,c):\mathbf{c}_{T}\rightarrow \mathbf{c}^{T}   \label{eq:morphism}
\end{align}
in the following way:
\begin{enumerate}
  \item To each colored basic tangle in Figure \ref{figure:mor}, associate a morphism accordingly:
        (a) ${\rm id}_{\mathbf{u}}$; (b) $R_{\mathbf{u},\mathbf{v}}$; (c) $R_{\mathbf{v},\mathbf{u}}^{-1}$; (d) $\iota_{\mathbf{u}}$; (e) $\epsilon_{\mathbf{u}}$.
  \item Call a framed tangle {\it elementary} if it is the horizontal composite of some basic tangles. For a colored elementary framed
        tangle $(T,c)$, let $\mathcal{F}(T,c)$ be the tensor product of the morphisms associated to the basic pieces.
  \item It is always possible to decompose a general $T$ as a vertical composite of elementary
        tangles; set $\mathcal{F}(T,c)$ to be the vertical composite of the morphisms associated to the elementary layers.
\end{enumerate}

It is know that $\mathcal{F}(T,c)$ is independent of the choice of the decomposition, and called the {\it Reshetikhin-Turaev invariant} of the colored framed tangle $(T,c)$.

\medskip

Now let $\mathcal{C}=\mathcal{E}(\Gamma)$.
Due to the coincidence between ${\rm DW}_{3,\Gamma}^0$ and ${\rm RT}_{\mathcal{E}(\Gamma)}$ \cite{Fr94-2}, the DW invariant of a link $L$ can be expressed in terms of RT invariants of $(L,c)$ for various colorings $c$.
\begin{thm} \label{thm:main}
For each $n$-component framed link $L=K_{1}\sqcup\cdots\sqcup K_n$,
\begin{align*}
Z(L)=\frac{1}{\#\Gamma}\sum\limits_{c\in{\rm col}(L)}\mathcal{F}(L,c)\cdot\otimes_{j=1}^{n}\chi_{c(j)}\in E^{\otimes n},
\end{align*}
where ${\rm col}(L)$ is the set of maps $c:\{1,\ldots,n\}\rightarrow \Lambda$.
\end{thm}
This result is somehow well-known to experts. Here are some explanations.
The Reshetikhin-Turaev theorem asserts that a MTC $\mathcal{C}$ determines a 1-2-3 extended TQFT ${\rm RT}_{\mathcal{C}}$ which takes $\mathcal{C}$ at $S^1$. In \cite{Fr94-1,Fr94-2} Freed used {\it finite integral} to construct an extended TQFT (extending the one presented in Section 2.1) whose value on the circle is computed to be $\mathcal{E}(\Gamma)$. Consequently, the 3-manifold invariants in ${\rm DW}_{3,\Gamma}^0$ are equal to those in ${\rm RT}_{\mathcal{E}(\Gamma)}$. Alternatively, one may refer to Theorem 4.4 of \cite{Ch12-2}, where many explicit details were given.

\section{The computing method for arborescent links}

Let $\mathcal{T}_{2}^{2}$ denote the set of tangles $T$ with $\#s(T)=\#t(T)=2$. Besides the vertical composition $\ast$, there is another operation on $\mathcal{T}_{2}^{2}$, namely, the rotation $r$: for any $T\in\mathcal{T}_{2}^{2}$, let $r(T)$ denote the new tangle obtained from rotating $T$ counterclockwise by $\pi/2$. Let $\mathcal{T}_{{\rm ar}}\subset \mathcal{T}_{2}^{2}$ denote the smallest subset containing $[1]$ and closed under $\ast$ and $r$.

Elements of $\mathcal{T}_{{\rm ar}}$ are called {\it arborescent tangles}, each of which can be written as a word in the alphabet $[1], \ast, r$. An {\it arborescent link} is the closure of some arborescent tangle.

Some special arborescent tangles and links are more interesting.

For $s_{1},\ldots,s_{k}\in\mathbb{Z}$, define the {\it continued fraction} $[[s_{1},\ldots,s_{k}]]$ inductively by
$$[[s_{1}]]=s_{1}, \qquad   [[s_{1},\ldots,s_{k}]]=s_{k}-1/[[s_{1},\ldots,s_{k-1}]].$$
Associated to each $p/q\in\mathbb{Q}$ (with $(p,q)=1$), there is a {\it rational tangle}
$$[p/q]=r(\cdots r([s_1])\ast\cdots)\ast[s_k],$$
where the $s_{i}$'s are chosen so that $[[s_{1},\ldots,s_{k}]]=p/q$; it is known that up to equivalence $[p/q]$ does not depend on the choices.

A {\it Montesinos tangle} is the composite tangle $r([p_1/q_1])\ast\cdots\ast r([p_m/q_m])$ for some rational numbers $p_{i}/q_{i}$,
and its closure is called a {\it Montesinos link}.

\begin{figure} [h]
  \centering 
  \includegraphics[width=11.5cm]{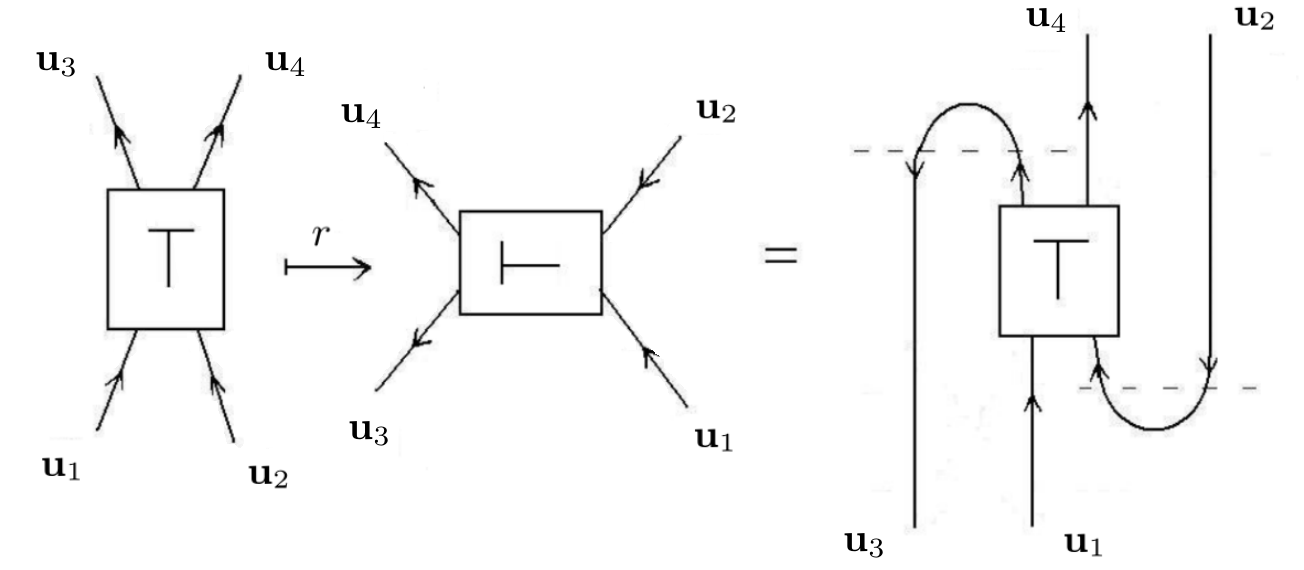}\\
  \caption{The rotation $r$.} \label{fig:rot}
\end{figure}

\begin{nota}
\rm Let $\mathcal{E}=\mathcal{E}(\Gamma)$. For $\mathbf{u},\mathbf{v}$, let $\mathcal{E}(\mathbf{u},\mathbf{v})$ denote the set of morphisms from $\mathbf{u}$ to $\mathbf{v}$.
\end{nota}

Let $(T,c)$ be a colored framed tangle, with $T\in\mathcal{T}_{2}^{2}$, as indicated on the left of Figure \ref{fig:rot}.
Associated to it is a morphism $\mathcal{F}(T,c)\in\mathcal{E}(\mathbf{u}_1\odot\mathbf{u}_2,\mathbf{u}_3\odot \mathbf{u}_4)$.
Figure \ref{fig:rot} shows that $\mathcal{F}(r(T,c))\in\mathcal{E}(\mathbf{u}_3^{\star}\odot \mathbf{u}_1,\mathbf{u}_4\odot \mathbf{u}_2^{\star})$
is equal to the composite
\begin{align*}
(\epsilon_{\mathbf{u}_3}\odot{\rm id}_{\mathbf{u}_4}\odot{\rm id}_{\mathbf{u}_2^{\star}})\circ({\rm id}_{\mathbf{u}_3^{\star}}\odot \mathcal{F}(T,c)\odot{\rm id}_{\mathbf{u}_2^{\star}})\circ({\rm id}_{\mathbf{u}_3^{\star}}\odot{\rm id}_{\mathbf{u}_1}\odot \iota_{\mathbf{u}_2}).
\end{align*}

For any $\mathbf{u}_{1}, \mathbf{u}_{2}, \mathbf{u}_{3}, \mathbf{u}_{4} \in\mathcal{E}$, define
\begin{align}
{\rm rot}={\rm rot}_{\mathbf{u}_1,\mathbf{u}_2}^{\mathbf{u}_3,\mathbf{u}_4}:\mathcal{E}(\mathbf{u}_1\odot \mathbf{u}_2,\mathbf{u}_3\odot\mathbf{u}_4)\rightarrow\mathcal{E}(\mathbf{u}_{3}^{\star}\odot\mathbf{u}_1,\mathbf{u}_4\odot\mathbf{u}_2^{\star})
\end{align}
by sending $F$ to the composite
$$(\epsilon_{\mathbf{u}_3}\odot{\rm id}_{\mathbf{u}_{4}}\odot{\rm id}_{\mathbf{u}_{2}^{\star}})\circ
({\rm id}_{\mathbf{u}_3^{\star}}\odot F\odot{\rm id}_{\mathbf{u}_2^{\star}})\circ({\rm id}_{\mathbf{u}_3^{\star}}\odot{\rm id}_{\mathbf{u}_1}\odot \iota_{\mathbf{u}_2}).$$
Then $\mathcal{F}(r(T,c))={\rm rot}(\mathcal{F}(T,c))$.

The functoriality and invariance of (\ref{eq:morphism}) immediately imply
\begin{thm} \label{thm:DW-of-arborescent}
For each colored arborescent tangle $(T,c)$ where $T$ is presented as a word $w([1],\ast,r)$, the morphism $\mathcal{F}(T,c)$ can be obtained from replacing $[1], \ast, r$ in $w$ by $R, \circ, {\rm rot}$, respectively.
\end{thm}

In particular, if $p/q=[[s_{1},\cdots,s_{k}]]$, then
\begin{align*}
\mathcal{F}([p/q],c)={\rm rot}(\cdots {\rm rot}(R^{s_1})\circ\cdots)\circ R^{s_{k}}.
\end{align*}

\begin{rmk}
\rm By $R$ we mean $R_{\mathbf{v},\mathbf{w}}$ for some $\mathbf{v},\mathbf{w}$ as determined by the coloring. Similarly for ${\rm rot}$.
\end{rmk}

For $\mathbf{v}\in\mathcal{E}$ and $f\in\mathcal{E}(\mathbf{v},\mathbf{v})$, let
\begin{align}
\overrightarrow{{\rm tr}}(f)=\sum\limits_{\mathbf{w}\in\Lambda}{\rm tr}\big(f_{\ast}:\mathcal{E}(\mathbf{w},\mathbf{v})\rightarrow\mathcal{E}(\mathbf{w},\mathbf{v})\big)\cdot\chi_{\mathbf{w}}\in E.
\end{align}
For $\alpha=\sum_{\mathbf{w}\in\Lambda}a_{\mathbf{w}}\chi_{\mathbf{w}}\in E$, put
\begin{align}
|\alpha|=\sum\limits_{\mathbf{w}\in\Lambda}\dim\mathbf{w}\cdot a_{\mathbf{w}}.
\end{align}

\begin{lem} \label{lem:tangle-closure}
Suppose $(T,c)$ is a colored framed tangle such that $s(T)$ can be identified with $t(T)$ and $\mathbf{c}_{T}=\mathbf{c}^{T}$, so that we can get a colored framed link $(\overline{T},\overline{c})$ by taking closure.
Then
\begin{align*}
\mathcal{F}(\overline{T},\overline{c})=\big|\overrightarrow{{\rm tr}}(\mathcal{F}(T,c))\big|.
\end{align*}
\end{lem}

\begin{proof}
By Lemma 2.6 of \cite{RT91}, $\mathcal{F}(\overline{T},\overline{c})={\rm tr}_{q}\big(\mathcal{F}(T,c)\big)$, where ${\rm tr}_{q}(\cdot)$,
in our setting and notation, is equal to $\big|\overrightarrow{{\rm tr}}(\cdot)\big|$, as can be checked.
\end{proof}

We now explain how to practically compute the RT invariant of a colored arborescent link. Suppose $T\in\mathcal{T}_{\rm ar}$ is presented as a word $w([1],\ast,r)$, and $c$ is a coloring of $T$ which induces a coloring $\overline{c}$ of $\overline{T}$.
Then by Theorem \ref{thm:DW-of-arborescent}, $\mathcal{F}(T,c)=w(R,\circ,{\rm rot})$.
\begin{enumerate}
  \item For all $\mathbf{u},\mathbf{v}\in\Lambda$, choose an isomorphism
        \begin{align}
        q_{\mathbf{u},\mathbf{v}}:\mathbf{u}\odot\mathbf{v}\cong
        \bigoplus\limits_{\mathbf{w}\in\Lambda}N_{\mathbf{u},\mathbf{v}}^{\mathbf{w}}\cdot \mathbf{w},
        \end{align}
        where $N_{\mathbf{u},\mathbf{v}}^{\mathbf{w}}\cdot\mathbf{w}$ stands for the direct sum of $N_{\mathbf{u},\mathbf{v}}^{\mathbf{w}}$ copies of $\mathbf{w}$.

        For $\mathbf{u}_{1},\mathbf{u}_{2},\mathbf{u}_{3},\mathbf{u}_{4}\in\Lambda$, let $\Phi=\Phi_{\mathbf{u}_{1},\mathbf{u}_{2}}^{\mathbf{u}_{3},\mathbf{u}_{4}}$ be the composite
        \begin{align}
        \mathcal{E}(\mathbf{u}_{1}\odot\mathbf{u}_{2},\mathbf{u}_{3}\odot \mathbf{u}_{4})&\cong\bigoplus\limits_{\mathbf{w}\in\Lambda}
        \mathcal{E}(N_{\mathbf{u}_{1},\mathbf{u}_{2}}^{\mathbf{w}}\cdot\mathbf{w},\ N_{\mathbf{u}_{3},\mathbf{u}_{4}}^{\mathbf{w}}\cdot\mathbf{w}) \label{eq:iso1} \\
        &\cong\bigoplus\limits_{\mathbf{w}\in\Lambda}
        \mathcal{M}(N_{\mathbf{u}_{1},\mathbf{u}_{2}}^{\mathbf{w}},N_{\mathbf{u}_{3},\mathbf{u}_{4}}^{\mathbf{w}}) \label{eq:iso2},
        \end{align}
        where (\ref{eq:iso1}) sends $f$ to $q_{\mathbf{u}_{3},\mathbf{u}_{4}}\circ f\circ q_{\mathbf{u}_{1},\mathbf{u}_{2}}^{-1}$, and (\ref{eq:iso2}) comes from the fact that each endomorphism of a simple object $\mathbf{w}$ is a multiple of ${\rm id}_{\mathbf{w}}$. We shall abbreviate
        $\Phi_{\mathbf{u}_{1},\mathbf{u}_{2}}^{\mathbf{u}_{3},\mathbf{u}_{4}}$ to $\Phi$ whenever there is no risk of confusion.
  \item Define ${\rm ROT}={\rm ROT}_{\mathbf{u}_{1},\mathbf{u}_{2}}^{\mathbf{u}_{3},\mathbf{u}_{4}}$ to be the unique map making
        the diagram
        \begin{align*}
        \xymatrix{
        \mathcal{E}(\mathbf{u}_{1}\odot\mathbf{u}_{2},\mathbf{u}_{3}\odot\mathbf{u}_{4})\ar[r]^{{\rm rot}} \ar[d]_{\Phi}
        &\mathcal{E}(\mathbf{u}_{3}^{\star}\odot\mathbf{u}_{1},\mathbf{u}_{4}\odot\mathbf{u}_{2}^{\star}) \ar[d]^{\Phi} \\
        \bigoplus\limits_{\mathbf{w}\in\Lambda}
        \mathcal{M}(N_{\mathbf{u}_{1},\mathbf{u}_{2}}^{\mathbf{w}},N_{\mathbf{u}_{3},\mathbf{u}_{4}}^{\mathbf{w}})\ar[r]^{{\rm ROT}} &\bigoplus\limits_{\mathbf{w}\in\Lambda}
        \mathcal{M}(N_{\mathbf{u}_{3}^{\star},\mathbf{u}_{1}}^{\mathbf{w}},N_{\mathbf{u}_{4},\mathbf{u}_{2}^{\star}}^{\mathbf{w}})
        }
        \end{align*}
        commutative.
        Choose base $\{u_{1,j}\}, \{u_{2,k}\}, \{u_{3,u}\}, \{u_{4,v}\}$ for $\mathbf{u}_1$, $\mathbf{u}_{2}$, $\mathbf{u}_{3}$, $\mathbf{u}_{4}$, respectively, and let $\{u_{2,k}^{\ast}\}, \{u_{3,u}^{\ast}\}$ be the dual base for $\mathbf{u}_{2}^\star, \mathbf{u}_{3}^\star$, respectively.

        Suppose $F:\mathbf{u}_{1}\odot \mathbf{u}_{2}\rightarrow \mathbf{u}_{3}\odot \mathbf{u}_{4}$ is expressed as
        \begin{align}
        F(u_{1,j}\otimes u_{2,k})=\sum\limits_{u,v}\mathcal{T}(F)_{j,k}^{u,v}\cdot u_{3,u}\otimes u_{4,v}.
        \end{align}
        Then ${\rm rot}(F)$ can be expressed as
        \begin{align*}
        u_{3,i}^{\ast}\otimes u_{1,j}&\mapsto
        \sum\limits_{k} u_{3,i}^{\ast}\otimes u_{1,j}\otimes u_{2,k}\otimes u_{2,k}^{\ast} \nonumber \\
        &\mapsto\sum\limits_{k,u,v}\mathcal{T}(F)_{j,k}^{u,v}\cdot u_{3,i}^{\ast}\otimes u_{3,u}\otimes u_{4,v}\otimes u_{2,k}^{\ast} \nonumber \\
        &\mapsto\sum\limits_{k,v}\mathcal{T}(F)_{j,k}^{i,v}\cdot u_{4,v}\otimes u_{2,k}^{\ast}.
        \end{align*}
        Hence with respect to the base $\{u_{3,i}^{\ast}\}, \{u_{1,j}\}, \{u_{4,k}\}, \{u_{2,v}^{\ast}\}$,
        $\textrm{rot}(F)$ is represented by the matrix $\mathcal{T}({\rm rot}(F))$, with
        \begin{align}
        \mathcal{T}({\rm rot}(F))_{i,j}^{k,v}=\mathcal{T}(F)_{j,v}^{i,k},  \label{eq:rot}
        \end{align}
        and the linear map ${\rm ROT}$ can be determined through this relation.
  \item For $A=\oplus_{\mathbf{w}\in\Lambda}A_{\mathbf{w}}\in\bigoplus_{\mathbf{w}\in\Lambda}\mathcal{M}(N^{\mathbf{w}},N^{\mathbf{w}})$, put
        \begin{align}
        \overrightarrow{{\rm Tr}}(A)=\sum\limits_{\mathbf{w}\in\Lambda}{\rm tr}(A_{\mathbf{w}})\cdot\chi_{\mathbf{w}}\in E.
        \end{align}
        Then by Lemma \ref{lem:tangle-closure},
        \begin{align}
        \mathcal{F}(\overline{T},\overline{c})=|\overrightarrow{{\rm tr}}(w(R,\circ,{\rm rot}))|
        =\big|\overrightarrow{{\rm Tr}}(w(\Phi(R),\circ,{\rm ROT}))\big|.
        \end{align}
\end{enumerate}


\section{$\Gamma=\mathbb{Z}/p\mathbb{Z}\rtimes\mathbb{Z}/(p-1)\mathbb{Z}$ for an odd prime $p$}

Let $\mathbb{Z}_p=\mathbb{Z}/p\mathbb{Z}$ and $\mathbb{Z}_{p-1}=\mathbb{Z}/(p-1)\mathbb{Z}$.

Fix a generator $r$ of $\mathbb{Z}_p^\times\cong\mathbb{Z}_{p-1}$. Let $\Gamma=\mathbb{Z}_p\rtimes\mathbb{Z}_{p-1}$ be determined by the homomorphism $\mathbb{Z}_{p-1}\to{\rm Aut}(\mathbb{Z}_{p})\cong\mathbb{Z}_{p}^\times$ sending $1$ to $r$. A presentation is
$$\Gamma=\langle\alpha,\beta\mid\alpha^p=\beta^{p-1}=1,\ \beta\alpha=\alpha^r\beta\rangle.$$
We often denote an element of $\Gamma$ by $\alpha^a\beta^b$ with $a\in\mathbb{Z}_p$, $b\in\mathbb{Z}_{p-1}$.

There are $p-1$ nontrivial conjugacy classes, namely,
\begin{align*}
{\rm Con}(\alpha)&=\{\beta^k\lrcorner\alpha=\alpha^{r^k}\colon k\in\mathbb{Z}_{p-1}\};  \\
{\rm Con}(\beta^v)&=\{\alpha^n\lrcorner\beta^v=\alpha^{(1-r^v)n}\beta^v\colon n\in\mathbb{Z}_{p}\}, \qquad  0\ne v\in\mathbb{Z}_{p-1}.
\end{align*}
The centralizers of the representatives are
$${\rm Cen}(\alpha)=\langle\alpha\rangle; \qquad {\rm Cen}(\beta^v)=\langle\beta\rangle, \ \  0\ne v\in\mathbb{Z}_{p-1}.$$

Let
$$\check{p}=\frac{p-1}{2}, \qquad  \zeta=\exp\Big(\frac{2\pi i}{p}\Big),  \qquad \xi=\exp\Big(\frac{2\pi i}{p-1}\Big).$$
It makes sense to write $\xi^k$ for $k\in\mathbb{Z}_{p-1}$ and $\zeta^s$ for $s\in\mathbb{Z}_p$.

\subsection{Simple objects}

Up to equivalence there is a unique $(p-1)$-dimensional irreducible representation
$\tilde{\mathbf{c}}:\Gamma\to{\rm GL}(p-1,\mathbb{C})$, which is determined by
\begin{align}
\alpha\mapsto\left[\begin{array}{ccc} \zeta^{r^{-1}} & \ & \  \\ \ & \ddots & \ \\ \ & \ & \zeta^{r^{-(p-1)}} \end{array}\right], \qquad
\beta\mapsto\left[\begin{array}{cccc} 0 & \cdots & 0 & 1 \\ 1 & \cdots & 0 & 0 \\
\vdots & \ddots & \vdots & \vdots \\ 0 & \cdots & 1 & 0 \end{array}\right],
\end{align}
so that $\beta(\mathfrak{e}_i)=\mathfrak{e}_{i+1}$, where $\mathfrak{e}_i$ is the $(p-1)$-dimensional column vector whose $i$-th entry is $1$ and the other entries are all $0$.
It is easy to find
$$\chi_{\mathbf{c}}(\alpha^x\beta^y)=(p\delta^0_x-1)\delta^0_y.$$
For each $k\in\mathbb{Z}_{p-1}$, there is a $1$-dimensional representation
\begin{align}
\tilde{\mathbf{d}}_k:\Gamma\to\mathbb{C}^\times, \qquad \alpha^u\beta^v\mapsto\xi^{kv}.
\end{align}
Let $\mathbf{c}$ (resp. $\mathbf{d}_k)$ denote the object of $\mathcal{E}(\Gamma)$ given by $\tilde{\mathbf{c}}$ (resp. $\tilde{\mathbf{d}}_k$) sitting at $e$.

For $j\in\mathbb{Z}_p$, set
\begin{align}
\mathbf{a}^j={\bigoplus}_{a=1}^{p-1}\mathbb{C}\langle\alpha^a\rangle,
\end{align}
by which we mean $\mathbf{a}^j_{\alpha^a}=\mathbb{C}$ for each $a$, and define the action by
$$\alpha^x\beta^y:\alpha^a\mapsto \zeta^{jx/(r^ya)}\cdot\alpha^{r^ya}.$$

For $0\ne v\in\mathbb{Z}_{p-1}$ and $\ell\in\mathbb{Z}_{p-1}$, set
\begin{align}
\mathbf{b}^\ell_v={\bigoplus}_{b=1}^p\mathbb{C}\langle\alpha^b\beta^v\rangle,
\end{align}
and define the action by
$$\alpha^x\beta^y:\alpha^b\beta^v\mapsto \xi^{\ell y}\cdot\alpha^{r^yb+x(1-r^v)}\beta^v.$$

Each simple object of $\mathcal{E}(\Gamma)$ is isomorphic to exactly one member in
$$\Lambda=\{\mathbf{c}\}\cup\{\mathbf{d}_k\colon k\in\mathbb{Z}_{p-1}\}\cup\{\mathbf{a}^j\colon j\in\mathbb{Z}_p\}\cup
\{\mathbf{b}_v^{\ell}\colon 0\ne v\in\mathbb{Z}_{p-1},\ \ell\in\mathbb{Z}_{p-1}\}.$$
By the definition (\ref{eq:def-dim}),
\begin{align}
\dim\mathbf{c}=p-1, \qquad \dim \mathbf{d}_k=1, \qquad \dim\mathbf{a}^j=p-1, \qquad \dim\mathbf{b}^\ell_v=p.  \label{eq:dim}
\end{align}
For the dual, we have
\begin{align}
\mathbf{c}^\star\cong\mathbf{c}, \qquad \mathbf{d}_k^\star\cong\mathbf{d}_{-k}, \qquad
(\mathbf{a}^j)^\star\cong \mathbf{a}^j, \qquad (\mathbf{b}^\ell_v)^\star\cong \mathbf{b}_{-v}^{-\ell}.
\end{align}
Precisely,
\begin{itemize}
  \item the isomorphism $\mathbf{c}\cong\mathbf{c}^\star$ is given by $\mathfrak{e}_i\mapsto\mathfrak{e}_{i+\check{p}}^\ast$;
  \item the isomorphism $\mathbf{d}_{-k}\cong\mathbf{d}_{k}^\star$ is obvious;
  \item the isomorphism $\mathbf{a}^j\cong(\mathbf{a}^j)^\star$ is given by
  \begin{align}
  \mathbf{a}^j_{\alpha^a}\to(\mathbf{a}^j)^\star_{\alpha^a}=(\mathbf{a}^j_{\alpha^{-a}})^\ast, \qquad \alpha^a\mapsto(\alpha^{-a})^\ast; \label{eq:dual-a}
  \end{align}
  \item the isomorphism $\mathbf{b}_{-v}^{-\ell}\cong(\mathbf{b}^\ell_v)^\star$ is given by
  \begin{align}
  (\mathbf{b}^{-\ell}_{-v})_{\alpha^b\beta^{-v}}\to ((\mathbf{b}^\ell_v)^\star)_{\alpha^b\beta^{-v}}
  =((\mathbf{b}^\ell_v)_{\alpha^{-r^v}\beta^v})^\ast, \qquad \alpha^b\beta^{-v}\mapsto(\alpha^{-r^vb}\beta^v)^\ast. \label{eq:dual-b}
  \end{align}
\end{itemize}

\subsection{$\odot$ and $R$}

We shall decompose $\mathbf{b}^{\ell_1}_{v_1}\odot \mathbf{b}^{\ell_2}_{v_2}$ as a direct sum of objects in $\Lambda$, and then determine $\Phi\big(R_{\mathbf{b}^{\ell_1}_{v_1},\mathbf{b}^{\ell_2}_{v_2}}\big)$ through this decomposition.

Let $\tilde{\ell}=\ell_1+\ell_2$ and $\tilde{v}=v_1+v_2$.

We do not deal with $\mathbf{u}\odot\mathbf{v}$ for $\mathbf{u}\ne\mathbf{b}_{v}^{\ell}$ or $\mathbf{v}\ne\mathbf{b}_{v}^{\ell}$, because a complete computation not only will occupy too many pages, but also is unnecessary (see Remark \ref{rmk:final}). The same reason is taken account for in the next subsection.

\subsubsection{$\tilde{v}\ne 0$}

In $\mathbf{b}^{\ell_1}_{v_1}\odot \mathbf{b}^{\ell_2}_{v_2}$, the component supported at $\beta^{\tilde{v}}$ is $\bigoplus_{a=1}^p\mathbf{x}_a$, with
$$\mathbf{x}_a=\mathbb{C}\langle \alpha^{-r^{v_1}a}\beta^{v_1}\rangle\otimes\mathbb{C}\langle\alpha^a\beta^{v_2}\rangle,$$
and the action $\beta:\mathbf{x}_a\to\mathbf{x}_{ra}$ is the multiplication by $\xi^{\tilde{\ell}}$, which in particular, fixes $\mathbf{x}_0$.
Hence
\begin{align}
\mathbf{b}^{\ell_1}_{v_1}\odot\mathbf{b}^{\ell_2}_{v_2}\cong\mathbf{b}^{\tilde{\ell}}_{\tilde{v}}\oplus
{\bigoplus}_{\ell=1}^{p-1}\mathbf{b}^{\ell}_{\tilde{v}},  \label{eq:tensor-1}
\end{align}
under which
\begin{align}
\beta^{v_1}\otimes\beta^{v_2}&\mapsto \beta^{\tilde{v}}\in \mathbf{b}^{\tilde{\ell}}_{\tilde{v}},  \label{eq:tensor-1-1} \\
\alpha^{-r^{v_1+u}}\beta^{v_1}\otimes\alpha^{r^u}\beta^{v_2}&\mapsto\oplus_{\ell=1}^{p-1}\xi^{(\ell-\tilde{\ell})u}\beta^{\tilde{v}}
\in{\bigoplus}_{\ell=1}^{p-1}\mathbf{b}^{\ell}_{\tilde{v}}. \label{eq:tensor-1-2}
\end{align}
Under the action of $R_{\mathbf{b}_{v_1}^{\ell_1},\mathbf{b}^{\ell_2}_{v_2}}$,
\begin{align*}
\beta^{v_1}\otimes\beta^{v_2}&\mapsto\xi^{\ell_2v_1}\cdot\beta^{v_2}\otimes\beta^{v_1}, \\
\alpha^{-r^{v_1+u}}\beta^{v_1}\otimes\alpha^{r^u}\beta^{v_2}
&\mapsto\xi^{\ell_2v_1}\cdot\alpha^{r^{u+\tilde{v}}}\beta^{v_2}\otimes\alpha^{-r^{v_1+u}}\beta^{v_1}.
\end{align*}
Hence through the decomposition (\ref{eq:tensor-1}),
\begin{align}
\Phi\big(R_{\mathbf{b}_{v_1}^{\ell_1},\mathbf{b}^{\ell_2}_{v_2}}\big)=\big(\xi^{\ell_2v_1}|\mathbf{b}^{\tilde{\ell}}_{\tilde{v}}\big)
\oplus\oplus_{\ell=1}^{p-1}\big(\xi^{(\ell-\tilde{\ell})\check{p}+(\ell-\ell_1)v_1}|\mathbf{b}^{\ell}_{\tilde{v}}\big),
\end{align}
where for a simple object $\mathbf{b}$, we use $(\lambda|\mathbf{b})$ to denote the morphism $\lambda\cdot{\rm id}_{\mathbf{b}}$.

\subsubsection{$v_1=v=-v_2\ne 0$}

In $\mathbf{b}^{\ell_1}_{v}\odot\mathbf{b}^{\ell_2}_{-v}$, the component supported at $\alpha^a$ is $\bigoplus_{b=1}^n\mathbf{y}_b$, with
$$\mathbf{y}_b=\mathbb{C}\langle\alpha^{a-r^vb}\beta^v\rangle\otimes\mathbb{C}\langle\alpha^b\beta^{-v}\rangle,$$
and the action $\alpha:\mathbf{y}_b\to \mathbf{y}_{b+(1-r^{-v})}$ is trivial.
Hence
\begin{align}
\mathbf{b}^{\ell_1}_{v}\odot \mathbf{b}^{\ell_2}_{-v}&\cong \mathbf{c}\oplus \mathbf{d}_{\tilde{\ell}}\oplus{\bigoplus}_{j=1}^{p}\mathbf{a}^{j}, \label{eq:tensor-2}
\end{align}
under which
\begin{align}
\alpha^{r^u-r^va}\beta^{v}\otimes\alpha^a\beta^{-v}&\mapsto \oplus_{j=1}^p\xi^{-\tilde{\ell}u}\zeta^{\frac{ja}{(1-r^{-v})r^u}}\alpha^{r^u}
\in{\bigoplus}_{j=1}^{p}\mathbf{a}^{j},  \label{eq:tensor-2-1} \\
\alpha^{-r^{v}a}\beta^{v}\otimes\alpha^{a}\beta^{-v}&\mapsto
\Big({\sum}_{i=1}^{p-1}\xi^{-\tilde{\ell}i}\zeta^{\frac{ar^{-i}}{1-r^{-v}}}\mathfrak{e}_i\Big)\oplus 1
\in\mathbf{c}\oplus \mathbf{d}_{\tilde{\ell}}.  \label{eq:tensor-2-2}
\end{align}
Under the action of $R_{\mathbf{b}_{v}^{\ell_1},\mathbf{b}^{\ell_2}_{-v}}$,
\begin{align*}
\alpha^{r^u-r^va}\beta^{v}\otimes\alpha^a\beta^{-v}&\mapsto\xi^{\ell_2v}\cdot\alpha^{a+r^u(1-r^{-v})}\beta^{-v}\otimes\alpha^{r^u-r^va}\beta^{v}, \\ \alpha^{-r^{v}a}\beta^{v}\otimes\alpha^{a}\beta^{-v}&\mapsto\xi^{\ell_2v}\cdot\alpha^{a}\beta^{-v}\otimes\alpha^{-r^{v}a}\beta^{v}.
\end{align*}
Hence through the decomposition (\ref{eq:tensor-2}),
\begin{align}
\Phi\big(R_{\mathbf{b}_{v}^{\ell_1},\mathbf{b}^{\ell_2}_{-v}}\big)=(\xi^{\ell_2v}|\mathbf{c})\oplus(\xi^{\ell_2v}|\mathbf{d}_{\tilde{\ell}})
\oplus\oplus_{j=1}^{p}\big(\xi^{\ell_2v}\zeta^{\frac{j}{1-r^v}}|\mathbf{a}^j\big).
\end{align}

\subsection{$\mathcal{T}(F)$}

We only provide information on $\mathcal{T}(F)$ for morphisms $F:\mathbf{b}^{\ell_1}_{v_1}\odot \mathbf{b}^{\ell_2}_{v_2}\to \mathbf{b}^{k_1}_{u_1}\odot \mathbf{b}^{k_2}_{u_2}$ with $(u_j,k_j)\in\{(v_j,\ell_j),(v_{3-j},\ell_{3-j})\}$.
Also, let $\tilde{\ell}=\ell_1+\ell_2$, $\tilde{v}=v_1+v_2$.

\subsubsection{$\tilde{v}\ne 0$}

Suppose $F:\mathbf{b}^{\ell_1}_{v_1}\odot \mathbf{b}^{\ell_2}_{v_2}\to \mathbf{b}^{k_1}_{u_1}\odot \mathbf{b}^{k_2}_{u_2}$ is a morphism with
$$\Phi(F)=\big(A|\mathbf{b}_{\tilde{v}}^{\tilde{\ell}}\oplus\mathbf{b}_{\tilde{v}}^{\tilde{\ell}}\big)\oplus
\widehat{\oplus}_{\ell=1}^{p-1}(b_\ell|\mathbf{b}_{\tilde{v}}^\ell),$$
where $A=[a_{ij}]_{2\times 2}$, and $\widehat{\oplus}$ means to take direct sum over $\ell\in\mathbb{Z}_{p-1}\setminus\{\tilde{\ell}\}$.

We shall express $\mathcal{T}(F)$ in terms of $\Phi(F)$ and vice versa.

Now $\beta^{v_1}\otimes\beta^{v_2}$ is sent by (\ref{eq:tensor-1-1}) to
$\beta^{\tilde{v}}$ in the first $\mathbf{b}_{\tilde{v}}^{\tilde{\ell}}$, and then sent by $\Phi(F)$ to $a_{11}\beta^{\tilde{v}}\oplus a_{12}\beta^{\tilde{v}}$;
if its image under the inverse of (\ref{eq:tensor-1}) is
$$a_{11}\alpha^{x(1-r^{u_1})}\beta^{u_1}\otimes\alpha^{x(1-r^{u_2})}\beta^{u_2}+
\sum\limits_{u=1}^{p-1}s_u\cdot\alpha^{x(1-r^{\tilde{v}})-r^{u_1+u}}\beta^{u_1}\otimes\alpha^{r^u}\beta^{u_2},$$
then by (\ref{eq:tensor-1-2}),
$$\sum\limits_{u=1}^{p-1}s_u\xi^{(\ell-\tilde{\ell})u}=\delta^{\ell}_{\tilde{\ell}}a_{12}.$$
This holds for all $\ell$ if and only if
\begin{align}
(p-1)s_u=a_{12}.
\end{align}

Moreover, for $b(1-r^{v_2})\ne 1-r^{\tilde{v}}$,
\begin{align*}
\alpha^{b-r^{v_1}}\beta^{v_1}\otimes\alpha\beta^{v_2}\stackrel{(\ref{eq:tensor-1-2})}\mapsto\ \oplus_{\ell=1}^{p-1}\alpha^{b}\beta^v
\stackrel{\Phi(F)}\mapsto a_{21}\alpha^{b}\beta^{v}\oplus\oplus_{\ell=1}^{p-1}b_\ell\alpha^{b}\beta^{\tilde{v}},
\end{align*}
with $b_{\tilde{\ell}}=a_{22}$; suppose it is sent by the inverse of (\ref{eq:tensor-1}) to
$$a_{21}\alpha^{b'(1-r^{u_1})}\beta^{u_1}\otimes\alpha^{b'(1-r^{u_2})}\beta^{u_2}+
\sum\limits_{w=1}^{p-1}c_{w}\cdot\alpha^{b-r^{u_1+w}}\beta^{u_1}\otimes\alpha^{r^{w}}\beta^{u_2},$$
with $b'(1-r^{\tilde{v}})=b$.
By (\ref{eq:tensor-1-2}),
\begin{align}
b_\ell=\sum\limits_{w=1}^{p-1}c_w\xi^{(\ell-\tilde{\ell})w}.
\end{align}
This holds for all $\ell$ if and only if
\begin{align*}
c_{w}=\frac{1}{p-1}\sum\limits_{\ell=1}^{p-1}\xi^{(\tilde{\ell}-\ell)w}b_{\ell}.
\end{align*}

Therefore,
\begin{align}
\mathcal{T}(F)_{x(1-r^{v_1}),x(1-r^{v_2})}^{x(1-r^{u_1}),x(1-r^{u_2})}&=a_{11}, \\
\mathcal{T}(F)_{x(1-r^{v_1}),x(1-r^{v_2})}^{x(1-r^{\tilde{v}})-r^{u_1+u},r^{u}}&=\frac{a_{12}}{p-1}, \\
\mathcal{T}(F)_{b'(1-r^{\tilde{v}})-r^{v_1+u},r^{u}}^{b'(1-r^{u_1}),b'(1-r^{u_2})}&=a_{21}, \\
\mathcal{T}(F)_{b-r^{v_1+u},r^{u}}^{b-r^{u_1+t},r^{t}}&=\frac{1}{p-1}\sum\limits_{\ell=1}^{p-1}\xi^{(\tilde{\ell}-\ell)t}b_{\ell}.
\end{align}
In the other direction (remembering that $a_{22}=b_{\tilde{\ell}}$),
\begin{align}
a_{11}=\mathcal{T}(F)^{0,0}_{0,0}, \qquad  a_{12}&=(p-1)\mathcal{T}(F)^{-r^{u_1},1}_{0,0},  \qquad a_{21}=\mathcal{T}(F)^{0,0}_{-r^{v_1},1}, \\
b_\ell&=\sum_{t=1}^{p-1}\xi^{(\ell-\tilde{\ell})t}\mathcal{T}(F)^{-r^{u_1+t},r^t}_{-r^{v_1},1}.
\end{align}

\subsubsection{$v_1=v=-v_2\ne0$}

Let $(u,k)=(v,\ell)$ or $(u,k)=(-v,-\ell)$. The computation is similar as above; note that $\tilde{\ell}=0$ in the present case.

Suppose $F:\mathbf{b}_{v}^{\ell}\odot \mathbf{b}_{-v}^{-\ell}\to \mathbf{b}_{u}^{k}\odot\mathbf{b}_{-u}^{-k}$ is a morphism with
$$\Phi(F)=(c|\mathbf{c})\oplus(d|\mathbf{d}_{0})\oplus\oplus_{j=1}^p(a_j|\mathbf{a}^j).$$

For $b\ne 0$,
\begin{align*}
\alpha^{b-r^va}\beta^v\otimes\alpha^{a}\beta^{-v}&\stackrel{(\ref{eq:tensor-2-1})}\mapsto\oplus_{j=1}^p\zeta^{\frac{ja}{(1-r^{-v})b}}\alpha^b
\stackrel{\Phi(F)}\mapsto\oplus_{j=1}^pa_j\zeta^{\frac{ja}{(1-r^{-v})b}}\alpha^b \\
&\ \mapsto\sum\limits_{x=1}^{p}\mu_{x}\cdot\alpha^{b-r^{u}x}\beta^{u}\otimes\alpha^{x}\beta^{-u},
\end{align*}
where the last map is the inverse of (\ref{eq:tensor-2-1}); from (\ref{eq:tensor-2-1}) we see
\begin{align}
\sum\limits_{x=1}^{p}\mu_{x}\zeta^{\frac{jx}{(1-r^{-u})b}}=a_j\zeta^{\frac{ja}{(1-r^{-v})b}}   \qquad \text{for\ all\ }j,
\end{align}
which is equivalent to
\begin{align*}
\mu_{x}=\frac{1}{p}\sum\limits_{j=1}^pa_j\zeta^{\frac{j}{b}(\frac{a}{1-r^{-v}}-\frac{x}{1-r^{-u}})}.
\end{align*}

For the remaining case,
\begin{align*}
\alpha^{-r^va}\beta^v\otimes\alpha^a\beta^{-v}&\stackrel{(\ref{eq:tensor-2-2})}\mapsto
\sum\limits_{i=1}^{p-1}\zeta^{\frac{ar^{-i}}{1-r^{-v}}}\mathfrak{e}_i\oplus 1
\stackrel{\Phi(F)}\mapsto c\sum\limits_{i=1}^{p-1}\zeta^{\frac{ar^{-i}}{1-r^{-v}}}\mathfrak{e}_i\oplus d \\
&\ \mapsto\sum\limits_{x=1}^{p}\eta_{x}\cdot\alpha^{-r^{u}x}\beta^{u}\otimes\alpha^{x}\beta^{-u},
\end{align*}
where the last is the inverse of (\ref{eq:tensor-2-2}), so
\begin{align*}
d=\sum\limits_{x=1}^{p}\eta_{x}, \qquad
c=\zeta^{\frac{ar^{-i}}{r^{-v}-1}}\sum\limits_{x=1}^{p}\eta_{x}\zeta^{\frac{xr^{-i}}{1-r^{-u}}}  \quad  \text{for\ all\ }i;
\end{align*}
equivalently,
\begin{align*}
\eta_{x}=\frac{1}{p}\Big(d+c\sum\limits_{i=1}^{p-1}\zeta^{r^{-i}(\frac{a}{1-r^{-v}}-\frac{x}{1-r^{-u}})}\Big)
=\begin{cases} c\delta^{x}_{a}+(d-c)/p, &u=v, \\ c\delta^{x}_{-r^va}+(d-c)/p, &u=-v.\end{cases}
\end{align*}

Therefore,
\begin{align}
\mathcal{T}(F)_{-r^va,a}^{-r^{u}x,x}&=\begin{cases} c\delta^{x}_{a}+(d-c)/p, &u=v, \\ c\delta^{x}_{-r^va}+(d-c)/p, &u=-v.\end{cases} \label{eq:Phi-to-T-2-1} \\
\mathcal{T}(F)_{b-r^va,a}^{b-r^{u}x,x}&=\frac{1}{p}\sum\limits_{j=1}^pa_j\zeta^{\frac{j}{b}(\frac{a}{1-r^{-v}}-\frac{x}{1-r^{-u}})}.
\label{eq:Phi-to-T-2-2}
\end{align}

In the other direction,
\begin{align}
c&=\mathcal{T}(F)^{0,0}_{0,0}-\mathcal{T}(F)^{-r^u,1}_{0,0}, \label{eq:T-to-Phi-2-1}  \\
d&=\mathcal{T}(F)^{0,0}_{0,0}+(p-1)\mathcal{T}(F)^{-r^u,1}_{0,0},   \label{eq:T-to-Phi-2-2}  \\
a_j&=\sum_{x=0}^{p-1}\zeta^{\frac{jx}{1-r^{-u}}}\mathcal{T}(F)^{1-r^ux,x}_{1,0}.   \label{eq:T-to-Phi-2-3}
\end{align}

\subsection{Example}

Let $m=2k+1$. Let $c$ be the coloring of $[m]$ whose values at the left- and right foot of $[m]$ are respectively $\mathbf{u}:=\mathbf{b}_v^\ell$ and $\mathbf{u}^\star\cong\mathbf{b}_{-v}^{-\ell}$. Let
$$F=\mathcal{F}([m],c)=(R_{\mathbf{u},\mathbf{u}^\star}\circ R_{\mathbf{u}^\star,\mathbf{u}})^k\circ R_{\mathbf{u},\mathbf{u}^\star} :\mathbf{b}_v^\ell\odot\mathbf{b}_{-v}^{-\ell}\to\mathbf{b}_{-v}^{-\ell}\odot \mathbf{b}_v^\ell.$$
Then
$$\Phi(F)=(\xi^{-m\ell v}|\mathbf{c})\oplus(\xi^{-m\ell v}|\mathbf{d}_0)\oplus\oplus_{j=1}^p\big(\xi^{-m\ell v}\zeta^{j(k+\frac{1}{1-r^v})}|\mathbf{a}^j\big).$$
By (\ref{eq:Phi-to-T-2-1}), (\ref{eq:Phi-to-T-2-2}),
\begin{align}
\mathcal{T}(F)^{-r^{-v}x,x}_{-r^va,a}&=\xi^{-m\ell v}\delta^x_{-r^va},  \\
\mathcal{T}(F)^{b-r^{-v}x,x}_{b-r^va,a}&=\frac{\xi^{-m\ell v}}{p}\sum_{j=1}^p\zeta^{j(k+\frac{1}{1-r^v})}\zeta^{\frac{j}{b}(\frac{a}{1-r^{-v}}-\frac{x}{1-r^v})}
=\xi^{-m\ell v}\delta^{x+r^va}_{(k+1-kr^v)b}.
\end{align}

We shall compute ${\rm rot}(F):\mathbf{b}_v^\ell\odot\mathbf{b}_v^\ell\to\mathbf{b}_v^\ell\odot\mathbf{b}_v^\ell$.
Be careful that (\ref{eq:rot}) holds with respect to the base $\{u_{3,i}^{\ast}\}, \{u_{1,j}\}, \{u_{4,k}\}, \{u_{2,v}^{\ast}\}$. Here we have used the dual (\ref{eq:dual-b}) to replace $(\mathbf{b}_{-v}^{-\ell})^\star$ by $\mathbf{b}_v^\ell$, so that the correct relation is
\begin{align}
\mathcal{T}({\rm rot}(F))^{j,a}_{s,t}=\mathcal{T}(F)^{-r^{-v}s,j}_{t,-r^{-v}a}.
\end{align}

Suppose $v\ne\check{p}$ so that $\tilde{v}=2v\ne 0$.
If
$$\Phi({\rm rot}(F))=(A'|\mathbf{b}_{2v}^{2\ell}\oplus\mathbf{b}_{2v}^{2\ell})\oplus\widehat{\oplus}_{s=1}^{p-1}(b'_s|\mathbf{b}_{2v}^s),$$
with $A'=[a'_{ij}]_{2\times 2}$, then
\begin{align*}
a'_{11}&=\mathcal{T}({\rm rot}(F))^{0,0}_{0,0}=\mathcal{T}(F)^{0,0}_{0,0}=\xi^{-m\ell v}, \\
a'_{12}&=(p-1)\mathcal{T}({\rm rot}(F))^{-r^v,1}_{0,0}=(p-1)\mathcal{T}(F)^{0,-r^v}_{0,-r^{-v}}=(p-1)\xi^{-m\ell v}\delta^{r^{-v}}_{1+1/k}, \\
a'_{21}&=\mathcal{T}({\rm rot}(F))^{0,0}_{-r^v,1}=\mathcal{T}(F)^{1,0}_{1,0}=\xi^{-m\ell v}\delta^{r^v}_{1+1/k}, \\
b'_s&=\sum_{t=1}^{p-1}\xi^{(s-2\ell)t}\mathcal{T}({\rm rot}(F))^{-r^{v+t},r^t}_{-r^v,1}
=\sum_{t=1}^{p-1}\xi^{(s-2\ell)t}\mathcal{T}(F)^{1,-r^{v+t}}_{1,-r^{t-v}} \\
&=\begin{cases} \xi^{(s-2\ell)z-m\ell v},&r^{\pm v}\ne 1+1/k, \\ 0,&\text{otherwise},\end{cases}
\qquad \text{with\ \ } r^z=\frac{kr^v-k-1}{(k+1)r^v-k}.
\end{align*}
The last line is computed as follows: If $(1,-r^{v+t},1,-r^{t-v})=(-r^{-v}x,x,-r^va,a)$ and $x=-r^va$, then $r^{v}=1$ which is impossible;
setting $(1,-r^{v+t},1,-r^{t-v})=(b-r^{-v}x,x,b-r^va,a)$ and $x+r^va=(k+1-kr^v)b$, we are led to
$$a=-r^{t-v}, \qquad b=1-r^t, \qquad x=-r^{t+v},  \qquad  ((k+1)r^v-k)r^t=kr^v-k-1.$$

Now suppose $v=\check{p}$, so that $\tilde{v}=0$ and $-r^{-v}=1$. Write
$$\Phi({\rm rot}(F))=(c'|\mathbf{c})\oplus(d'|\mathbf{d}_{2\ell})\oplus\oplus_{j=1}^p(a'_j|\mathbf{a}^j).$$
Noting that
$$\mathcal{T}({\rm rot}(F))^{0,0}_{0,0}=\mathcal{T}(F)^{0,0}_{0,0}=(-1)^\ell, \qquad
\mathcal{T}({\rm rot}(F))^{1,1}_{0,0}=\mathcal{T}(F)^{0,1}_{0,1}=(-1)^\ell\delta_{p\mid m},$$
we obtain
\begin{align*}
c'=(-1)^\ell(1-\delta_{p\mid m}),  \qquad    d'=(-1)^\ell p^{\delta_{p\mid m}}.
\end{align*}
Since $\mathcal{T}(F)^{1,1+y}_{0,y}=(-1)^\ell\delta^1_{-my}$, we have
\begin{align*}
a'_j=\sum_{y=1}^p\zeta^{\frac{jy}{2}}\mathcal{T}({\rm rot}(F))^{1+y,y}_{1,0}
=\sum_{y=1}^p\zeta^{\frac{jy}{2}}\mathcal{T}(F)^{1,1+y}_{0,y}
=\begin{cases} 0,&p\mid m, \\ (-1)^\ell\zeta^{-\frac{j}{2m}},&p\nmid m. \end{cases}
\end{align*}

\begin{figure}[h]
  \centering
  \includegraphics[width=5cm]{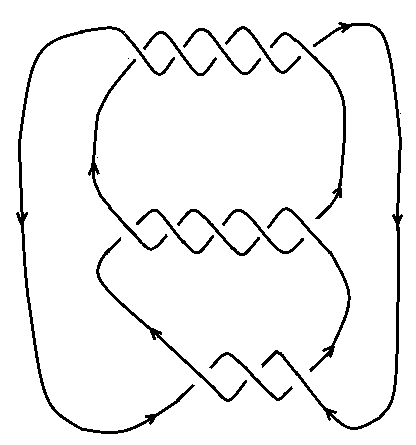}\\
  \caption{The odd classical pretzel knot $P(2k_1+1,2k_2+1,2k_3+1)$}\label{fig:pretzel}
\end{figure}

\bigskip

Let $K$ be the odd classical pretzel knot $P(n_1,n_2,n_3)$, with $n_i=2k_i+1$, as shown in Figure \ref{fig:pretzel};
let $n=n_1+n_2+n_3$. In the notation of Section 3, $K=\overline{T}$, with $T=r([n_1])\ast r([n_2])\ast r([n_3])$. Let $c$ denote the coloring of $T$ which induces the coloring $\mathbf{b}_{v}^\ell$ on $K$.

Suppose the coloring is $c=\mathbf{b}_{v}^\ell$ with $v\ne\check{p}$. We write
$$\big(\xi^{-n\ell v}B|\mathbf{b}_{2v}^{2\ell}\oplus\mathbf{b}_{2v}^{2\ell}\big)\oplus\widehat{\oplus}_{s=1}^{p-1}
        \big(\xi^{-n\ell v}b_s|\mathbf{b}_{2v}^s\big).$$
\begin{enumerate}
  \item If $r^{\pm v}\ne 1+1/k_i$ for $i=1,2,3$, then letting
        $$r^{w}=\prod_{i=1}^3\frac{k_ir^v-k_i-1}{(k_i+1)r^v-k_i},$$
        we have $B=I_2$, and $b_s=\xi^{(s-2\ell)w}$.
        Hence
        $$\overrightarrow{{\rm tr}}(\mathcal{F}(T,c))
        =\xi^{-n\ell v}\Big(\chi_{\mathbf{b}_{2v}^{2\ell}}+{\sum}_{s=1}^{p-1}\xi^{(s-2\ell)w}\chi_{\mathbf{b}_{2v}^{s}}\Big).$$
        By Lemma \ref{lem:tangle-closure},
        $$\mathcal{F}(K,\mathbf{b}_v^\ell)=p\xi^{-n\ell v}\big(1+\xi^{-2\ell w}(p-1)\delta_{p-1\mid w}\big)
        =\xi^{-n\ell v}p^{1+\delta^1_{r^w}}.$$
  \item If $r^v=1+1/k_i$ for some $i$ and $r^{-v}\ne 1+1/k_j$ for all $j$, then
        $B=\left[\begin{array}{cc} 1 & 0 \\ 0 & 0 \end{array}\right]$, and $b_s=0$ for all $s$.
        Hence $\overrightarrow{{\rm tr}}(\mathcal{F}(T,c))=\xi^{-n\ell v}\chi_{\mathbf{b}_{2v}^{2\ell}}$,
        so that $\mathcal{F}(K,\mathbf{b}_v^\ell)=p\xi^{-n\ell v}$.
  \item If $r^{-v}=1+1/k_i$ for some $i$ and $r^{v}\ne 1+1/k_j$ for all $j$, then similarly, $\mathcal{F}(K,\mathbf{b}_v^\ell)=p\xi^{-n\ell v}$.
  \item If $r^v=1+1/k_i$ for some $i$ and $r^{-v}=1+1/k_j$ for some $j$, then
        ${\rm tr}(B)=p^2$, and $b_s=0$. Hence $\mathcal{F}(K,\mathbf{b}_v^\ell)=p^2\xi^{-n\ell v}$.
\end{enumerate}

The Alexander polynomial of $K$ is (see \cite{Ch21} Example 4.3)
\begin{align}
\Delta_K(t)=1+(1+k_1+k_2+k_3+k_1k_2+k_1k_3+k_2k_3)(t+t^{-1}-2).
\end{align}
In case 1, $r^w=1$ is equivalent to $\Delta_K(r^v)=0\in\mathbb{Z}_p$. When $1+1/k_i\in\{r^{\pm v}\}$ for some $i$, $\Delta_K(r^v)\ne 0$ in case 2 and 3, and $\Delta_K(r^v)=0$ in case 4.

Consider the remaining case $v=\check{p}$.
\begin{itemize}
  \item If $p\nmid n_1n_2n_3$, then
        $$((-1)^\ell|\mathbf{c})\oplus((-1)^\ell|\mathbf{d}_{2\ell})\oplus\oplus_{j=1}^p
        \big((-1)^\ell\zeta^{-\frac{n_1n_2+n_1n_3+n_2n_3}{2n_1n_2n_3}j}\big|\mathbf{a}^j\big).$$
        Using $n_1n_2+n_1n_3+n_2n_3=-\Delta_K(-1)$, we obtain
        $$\overrightarrow{{\rm tr}}(\mathcal{F}(T,c))=(-1)^{\ell}\Big(\chi_{\mathbf{c}}+\chi_{\mathbf{d}_{2\ell}}
        +{\sum}_{j=1}^{p}\zeta^{\frac{\Delta_K(-1)}{2n_1n_2n_3}j}\chi_{\mathbf{a}^j}\Big),$$
        which results in
        $$\mathcal{F}(K,\mathbf{b}_{\check{p}}^\ell)=p(-1)^\ell\big(1+(p-1)\delta_{p\mid \Delta_K(-1)}\big)
        =(-1)^\ell p^{1+\delta_{p\mid\Delta_K(-1)}}.$$
  \item If $v=\check{p}$ and $\iota:=\{i\colon p\mid n_i\}>0$, then
        $$(0|\mathbf{c})\oplus((-1)^\ell p^\iota|\mathbf{d}_{2\ell})\oplus\oplus_{j=1}^p(0|\mathbf{a}^j),$$
        Hence $\overrightarrow{{\rm tr}}(\mathcal{F}(T,c))=(-1)^{\ell}p^\iota\chi_{\mathbf{d}_{2\ell}}$,
        implying $\mathcal{F}(K,\mathbf{b}_{\check{p}}^\ell)=(-1)^\ell p^\iota$. Note that $\iota=1+\delta_{p\mid\Delta_K(-1)}$ if $\iota\in\{1,2\}$, and $\iota=2+\delta_{p\mid\Delta_K(-1)}$ if $\iota=3$.
\end{itemize}

Put
\begin{align}
\epsilon(v)=\begin{cases} 1,&\Delta_K(r^v)=0\in\mathbb{Z}_p, \\ 0,&\text{otherwise};\end{cases} \qquad
\epsilon'=\begin{cases} 1,&p\mid n_1,n_2,n_3, \\ 0,&\text{otherwise}.\end{cases}
\end{align}
All the above cases can be uniformed as
$$\mathcal{F}(K,\mathbf{b}_v^\ell)=\xi^{-n\ell v}p^{1+\epsilon(v)+\epsilon'}.$$

Since $\chi_{\mathbf{b}_{v}^\ell}(\beta^v,\beta^y)=\xi^{\ell y}$, by Theorem \ref{thm:main} we have
\begin{align}
Z(K)(\beta^v,\beta^y)=\frac{1}{p(p-1)}\sum_{\ell=1}^{p-1}F(K,\mathbf{b}_{v}^\ell)\xi^{\ell y}=p^{\epsilon(v)+\epsilon'}\delta^{y}_{nv}.
\end{align}
This is the number of homomorphisms $\pi(K)\to\Gamma$ which take $\beta^v,\beta^y$ respectively at the meridian and the longitude.

\begin{rmk}\label{rmk:final}
\rm Taking conjugacy, one can easily obtain $Z(K)(\alpha^a\beta^v,\alpha^b\beta^y)$.

In contrast, $Z(K)(\alpha^a,\alpha^b)$ is not interesting,
since any homomorphism $\pi(K)\to\Gamma$ whose value at the meridian lies in $\langle\alpha\rangle$ must factor as $\pi(K)\to\langle\alpha\rangle\hookrightarrow\Gamma$.
\end{rmk}

\ \\
{\bf Conflict of interest}.
The author declares to have no commercial or associative interest that represents a conflict of interest in connection with the work submitted.

\ \\
Haimiao Chen \ \ \  \emph{chenhm@math.pku.edu.cn} \\
Department of Mathematics, Beijing Technology and Business University, Beijing, China.


\begin{thebibliography}{}


\bibitem{Ba01}
B. Bakalov,
\textsl{Lectures on tensor categories and modular functors}.
University Lecture Series, vol. 21. Providence, RI: American Mathematical Society; 2001.


\bibitem{BF20}
M. Boileau and S. Friedl,
The profinite completion of 3-manifold groups, fiberedness and the Thurston norm.
In: \textsl{What's next?--the mathematical legacy of William P. Thurston},
Ann. of Math. Stud., vol. 205. Princeton, NJ: Princeton Univ. Press; 2020. p. 21--44.


\bibitem{Ch12-1}
H.-M. Chen,
Applying TQFT to count regular coverings of Seifert 3-manifolds.
\textsl{J. Geom. Phys.} 62(6): 1347--1357, 2012.

\bibitem{Ch12-2}
H.-M. Chen,
The untwisted Dijkgraaf-Witten invariant of links. \\
arXiv:1209.4283v4.

\bibitem{Ch18}
H.-M. Chen,
On the groups of periodic links. 
arXiv:1805.02219.


\bibitem{Ch21}
H.-M. Chen,
Computing twisted Alexander polynomials for Montesinos links,
\textsl{Indian J. Pure. Appl. Math.} 52(2): 584--598, 2021.


\bibitem{DW90}
R. Dijkgraaf and E. Witten,
Topological gauge theories and group cohomology.
\textsl{Commun. Math. Phys.} 129(2): 393--429, 1990.


\bibitem{Ei07}
M. Eisermann,
Knot colouring polynomials.
\textsl{Pacific J. Math.} 231(2): 305--336, 2007.


\bibitem{Fe93}
K. Ferguson,
Link invariants associated to TQFT's with finite gauge groups.
\textsl{J. Knot Theory Ramifications} 2(1): 11--36, 1993.

\bibitem{FL21}
J. Flake and R. Laugwitz,
On the monoidal center of Deligne's category $\underline{\rm Rep}(S_t)$.
\textsl{J. Lond. Math. Soc. (2)} 103: 1153--1185, 2021.


\bibitem{FQ93}
D.S. Freed and F. Quinn,
Chern-Simons theory with finite gauge group.
\textsl{Commun. Math. Phys.} 156(3): 435--472, 1993.

\bibitem{Fr94-1}
D.S. Freed,
Higher algebraic structures and quantization.
\textsl{Commun. Math. Phys.} 159(2): 343--398, 1994.

\bibitem{Fr94-2}
D.S. Freed,
Quantum groups from path integrals.
In: \textsl{Particles and fields}, CRM Ser. Math. Phys. New York: Springer; 1994. p. 63--107.

\bibitem{Ha36}
P. Hall,
The Eulerian functions of a group.
\textsl{Quart. J. Math.} 7: 134--151, 1936.

\bibitem{Jo95}
G. Jones,
Enumeration of homomorphisms and surface-coverings.
\textsl{Quart. J. Math.} Oxford Ser. (2) 46 (184): 485--507, 1995.


\bibitem{Ko13}
T. Koberda,
Residual properties of fibered and hyperbolic 3-manifolds.
\textsl{Topol. Appl.} 160: 875--886, 2013.

\bibitem{LM00}
V. Liskovets and A. Mednykh,
Enumeration of subgroups in the fundamental groups of orientable circle bundles over surfaces.
\textsl{Commun. Algebra} 28(4): 1717--1738, 2000.


\bibitem{LR10}
D.D. Long and A.W. Reid,
Grothendieck's problem for 3-manifold groups.
\textsl{Groups Geom. Dynam.} 5(2): 479--499, 2010.


\bibitem{Lu95}
A. Lubotzky,
Counting finite index subgroups.
In: \textsl{Groups '93, Galway/St.Andrews}, vol. 2, London Mathametical Society, Lecture Note Series, vol.212.
Cambridge: Cambridge University Press; 1995. p. 368--404.


\bibitem{MS02}
D. Matei and A.I. Suciu,
Hall invariants, homology of subgroups, and characteristic varieties.
\textsl{Int. Math. Res. Not.} 2002(9): 465--503, 2002.

\bibitem{Mo15}
J.C. Morton,
Cohomological twisting of 2-linearization and extended TQFT.
\textsl{J. Homotopy Relat. Struct.} 10(2): 127--187, 2015.


\bibitem{PP06}
E. Pervova and C. Petronio,
On the existence of branced coverings between surfaces with prescribed branch data. I.
\textsl{Algebr. Geom. Topol.} 6: 1957--1985, 2006.


\bibitem{RT91}
N. Reshetikhin and V.G. Turaev,
Invariants of 3-manifolds via link polynomials and quantum groups.
\textsl{Invent. Math.} 103(3): 547--597, 1993.


\bibitem{Su01}
A.I. Suciu,
Fundamental groups of line arrangements: Enumerative aspects.
In: \textsl{Advances in algebraic geometry motivated by physics} (Lowell, MA, 2000), Contemp. Math. vol. 276.
Providence, RI: Amer. Math. Soc.; 2001. p. 43--79.

\bibitem{Tu09}
V.G. Turaev,
On certain enumeration problems in two-dimensional topology.
\textsl{Math. Res. Lett.} 16(3): 515--529, 2009.


\bibitem{Wa92}
M. Wakui,
On Dijkgraaf-Witten invariant for 3-manifolds.
\textsl{Osaka J. Math.} 29(4): 675--696, 1992.


\bibitem{Wi17}
G. Wilkes,
Virtual pro-$p$ properties of 3-manifold groups.
\textsl{J. Group Theory} 20(5): 999--1023, 2017.


\bibitem{WZ17}
H. Wilton and P. Zalesskii,
Distinguishing geometries using finite quotients.
\textsl{Geom. Topol.} 21(1): 345--384, 2017.


\bibitem{Ye92}
D.N. Yetter,
Topological quantum field theories associated to finite groups and crossed $G$-sets.
\textsl{J. Knot Theory Ramifications} 1(1): 1--20, 1992.

\end{thebibliography}
\end{document}